\magnification=1200

                       \bigskip
      \centerline {\bf  Odd Invariant Semidensity}
      \centerline {\bf and}
     \centerline { {\bf Divergence--like Operators on an
              Odd Symplectic Superspace}\footnote{$^*$}
                      {The work was supported in part by Suisse
          National Foundation and by INTAS--RFBR Grant No 95--0829}
                            }
                       \bigskip
                       \bigskip
       \centerline   {\bf O.M. Khudaverdian}
  \centerline   {\it  Laboratory of Computing Technique and Automation.
    Joint Institute for Nuclear Research}
    \centerline {{\it Dubna, Moscow Region  141980, Russia}
             \footnote{$^{**}$}
       {on leave of absence from Department of
     Theoretical Physics of Yerevan State
     University  375049  Yerevan, Armenia}
                      }
\centerline  {e-mails: "khudian@main1.jinr.dubna.su",
                          "khudian@sc2a.unige.ch"}
                      \bigskip
                       $$ $$
      {\it The divergence--like operator on an odd symplectic
   superspace which acts invariantly on a specially chosen
   odd vector field  is considered. This operator is used
  to construct an odd invariant semidensity
 in a geometrically clear way.
 The formula for this semidensity is similar to
  the formula of the mean curvature of hypersurfaces in
             Euclidean space.}
 \vfill
 \eject
     \medskip
    In this paper we construct a differential first order
     divergence-like operator on a superspace endowed
  with an odd symplectic structure.
  This operator is applied for constructing  invariant
  differential objects in this superspace.
   In particular  an odd semidensity
   invariant under the transformations preserving
 the  odd  symplectic structure of the superspace
   and the volume form is constructed in a geometrically clear way.

   The odd symplectic structure
   plays essential role in Lagrangian
    formulation of the BRST formalism
    (Batalin-Vilkovisky formalism ) [1].
    In supermathematics it is a natural counterpart of even symplectic
    structure [2,3]. On the other hand it has some
    "odd" features which have no analogues
    in usual mathematics. Canonical transformations
    preserving odd symplectic structure (non-degenerate
    odd Poisson bracket) do not preserve any volume form.
    This fact, indeed, is the reason why in the case of odd
      symplectic structure
       some invariant geometrical objects
      have no natural homologues
      as it is for the even symplectic structure,
       for which the superconstructions are
       rather straightforward generalizations
       of those for the symplectic structures in usual
       spaces. We consider some examples.

     In order to construct geometrical integration objects
    one needs to consider a pair, the volume form and the
    odd canonical structure.
    This pair is in fact the geometrical background for
    the formulation of Batalin-Vilkovisky formalism
    (See [4,5,6]). The so called $\Delta$ operator which plays
    essential role in this formalism can be defined in
    the following way: its action on the function $f$ is equal to the
    divergence with respect to the given volume form (defined by the
    action of the theory) of the Hamiltonian vector field,
     corresponding to the function $f$. It is a second order
    differential operator. In the case of even symplectic structure
    there exists the volume form which naturally corresponds
    to the canonical structure, and
    $\Delta$ operator is evidently equal to zero (Liouville Theorem).
    Even if the volume form is arbitrary,
     one arrives at a first order differential operator [4].

       The second example is the invariant
       volume form (density) which can be defined
       on the Lagrangian supersurfaces  embedded
       in the superspace endowed with an odd symplectic structure
     and a volume form. It is nothing
       but the integrand for the partition
       function in the space of fields and antifields
       in the BV formalism [5,6].

        We focus on the third example, on the problem of finding
        an analogue of Poincare-Cartan integral invariants
        for the odd symplectic structure.
        In usual mathematics to
         the Poincare-Cartan integral invariant
        (invariant volume form on the embedded surface)
        there corresponds the wedged power of a differential
        two-form which defines the symplectic structure.
        In case of even symplectic structure in spite of the
        fact that differential form has nothing in common
        with invariant integration objects, the superdeterminant
        of the two-form induced on the embedded supersurface
       from the two-form which defines the symplectic structure
	leads us to the Poincare--Cartan invariant [7].
        For the odd symplectic structure the situation
        is essentially different.
         In [8] there was considered the problem of constructing
    of the invariant densities for the superspace endowed with an
    odd symplectic structure.

     In the case of symplectic structure in usual mathematics
  as well as in the case of even symplectic structure
invariant densities are exhausted by
 densities \footnote{$^*$}{ Density
    is the object which defines a volume form
   on  embedded (super)surfaces.--- If
      $z=z(\zeta)$ is the local parametrization of a (super)surface
  then a density $A$
   is a function of
 $z(\zeta)$, ${\partial z\over\partial\zeta}$,
   ${\partial^2 z\over\partial\zeta\partial\zeta\dots}.$
   $A=A\left(z(\zeta),{\partial z\over\partial\zeta},
 {\partial^2\over\partial\zeta\partial\zeta},\dots,
   {\partial z^k\over\partial\zeta\dots\partial\zeta}\right)$
 subject to the condition that
  under  reparametrization
  $\zeta\rightarrow\zeta({\tilde\zeta}),\,A\rightarrow
A\cdot sdet\left({\partial\zeta\over\partial{\tilde\zeta}}
 \right)$.
   $A$ is a density of the rank $k$ if it depends on the
 tangent vectors of the order $\leq k$
 (i.e. on the derivatives of the order $\leq k$).
  In usual mathematics densities are the natural
 generalization of differential forms (if $k=1$). In supermathematics
 even in the case of $k=1$ this object is of more importance, since
 differential forms in supermathematics are no more integration
 objects [9,10,11].}  of the rank $k=1$
 which correspond to Poincare--Cartan integral invariants.
  This is not the case for the odd symplectic structure.
    One can show that on $(p.p)$-dimensional
   non-degenerate supersurfaces embedded in a superspace $E^{n.n}$
   which is endowed with a volume form and an odd
symplectic structure there are no invariant densities
   of the rank $k=1$ (except of the volume form itself),  and
        in the class of densities of the rank
       $k=2$ and of the weight $\sigma$
     (i.e. which are multiplied by the $\sigma$-th
  power of the superdeterminant of the reparametrization)
        there exists a unique (up to multiplication by a constant)
     semidensity ($\sigma={1\over 2}$)
      in the case of $p=n-1$ [8].
      In fact, in [8] this semidensity was
     constructed in a non-explicit way in  terms
    of dual densities: If $(n-1.n-1)$ supersurface
   $M^{n-1.n-1}$ is given by the equations:
      $f=0,\varphi=0$
   where $f$ is even function and $\varphi$ an odd function
then to this semidensity there corresponds the function
                   $$
                 {\tilde A}\Big\vert_{f=\varphi=0}=
           {1\over\sqrt{\{f,\varphi\}}}
                  \left(
                  \Delta f-
          {\{f,f\}\over 2\{f,\varphi\}}
                \Delta\varphi-
              {\{f,\{f,\varphi\}\}\over\{f,\varphi\}}-
                  {\{f,f\}\over 2\{f,\varphi\}^2}
              \{\varphi,\{f,\varphi\}\}
                       \right)\,.
                                       \eqno(1.1)
                        $$
   which depends on the second derivatives ($k=2$),
 is invariant under the transformations
preserving the odd symplectic structure and the volume form,
 and is multiplied by the square root of the corresponding
  Berezinian (superdeterminant) under the
   transformation
     $f\rightarrow$ $a f+\alpha\varphi$,
     $\varphi\rightarrow$ $\beta f+b\varphi$,
 (which does not change the supersurface $M$).
The semidensity (1.1) takes odd values.
 It is an exotic analogue of Poincare--Cartan
invariant: ${\tilde A}^2=0$, so it cannot be  integrated
 over supersurfaces.

    \smallskip
  To clarify the geometrical meaning of density (1.1),
  in this paper a special geometrical object
 for odd symplectic superspace is considered.

 As it was mentioned above, to a symplectic structure in usual space
 there corresponds a volume form (Liouville form):
                     $$
          d{\bf v}=\rho(x)dx^1...dx^{2n}=
            \sqrt{det\,(\Omega_{ik})}
                 dx^1...dx^{2n}
                                          \eqno (1.2)
                     $$
 where $\Omega=\Omega_{ik}dx^i\wedge dx^k$
 is  the closed non-degenerated two-form which defines
 the symplectic structure.

\noindent  The volume form (1.2) is preserved
 under canonical transformations
(i.e., the transformations preserving the two-form $\Omega$).
 If ${\bf X}$ $={X^i{\partial\over\partial x^i}}$ is
 an arbitrary vector field
one can consider its divergence:
                    $$
     div {\bf X}= {{\cal L}_{\bf X}d{\bf v}\over d{\bf v}}=
      {\partial X^i\over\partial x^i}+
      X^i {\partial\log\rho\over\partial x^i}\,.
                                                           \eqno (1.3)
                   $$
  In a symplectic space the canonical transformations
preserve not only  the volume form (1.2) but also its projection on
an arbitrary symplectic plane. Moreover, if
  ${\hat L}(z)$ is a projector-valued function such that
  ${\bf Im}{\hat L}(z)$ is a symplectic subspace (plane)
in $T_z E$ then it is easy to see that
to ${\hat L}(z)$ one can associate divergence-like invariant
 operator whose action on
an arbitrary vector field
 ${\bf X}$ is given by the expression
                        $$
 \partial({\hat L},{\bf X})=\left(
            {\partial X^k\over\partial x^i}
               +{1\over 2}
          X^m{\partial \Omega_{ip}\over\partial x^m}
                      \Omega^{pk}\right)
                      { L}^i_k\,.
                                                        \eqno (1.4)
                          $$
(Compare with (1.3) in case of ${\hat L}={\bf id}$).

 The formulae (1.2)--(1.4) have the straightforward generalization
to the case of even symplectic structure in superspace
(by changing determinant to superdeterminant and adding the
powers of $(-1)$ wherever necessary),
but it is not the case for the odd symplectic structure.

 Nevertheless, it turns out that the analogue of the formula
 (1.4) can be considered
 for odd symplectic structure in a case where ${\hat L}$
is a projector on $(1.1)$-dimensional symplectic subspace and
 ${\bf X}$ is the odd vector field which belongs to this subspace and
is symplectoorthogonal to itself. In the next two sections
we perform the corresponding constructions which are essentially
founded on the following remark.
 Let $E^{1.1}$ be $(1.1)$-dimensional odd symplectic
superspace and $(x,\theta)$ be Darboux coordinates on it:
 $\{x,\theta\}=1$,
  $\{x,x\}=0$ where $\{\,,\,\}$ is the odd
Poisson bracket (Buttin bracket) corresponding to
  this symplectic structure.
  Let ${\bf \Psi}$ be an odd vector field
in this superspace which is equal to
   $\Psi(x,\theta)\partial_\theta$
 in these Darboux coordinates.
 Then this vector field  has
the same form  and
 its divergence ($\partial\Psi(x,\theta)/\partial \theta$)
 remains the same in arbitrary Darboux coordinates.
(See Example 2).

 In the 4-th section we  consider a $(n-1.n-1)$-dimensional
  supersurface
embedded in an odd symplectic superspace $E^{n.n}$ endowed
with a volume form.
    The differential operator described above can be naturally
    applied to the  odd vector field which is defined only
 in the points
of this supersurface and is symplectoorthogonal
to itself and to this supersurface.
 In spite of the fact that this field is not defined on the
 whole superspace one can define the invariant "truncated divergence"
 of this vector field. The analogue of this construction for usual
 symplectic structure is trivial. On the other hand
this construction
can be considered as an analogue of the corresponding operator
acting on the vector field which is defined in a
 Euclidean (Riemannian)
space on the points of embedded hypersurface and which takes
values in the normal bundle to this surface.
    But the essential difference is that
the group of transformations which preserve metrics in Euclidean
 (Riemannian) space is
 exhausted by linear transformations
(the linear part of transformation defines uniquely
  all higher terms) and this is not true for symplectic case
where the group of canonical
 transformations is infinite-dimensional.

 In the 5-th Section we apply our geometrical construction to
obtain the formula for odd semidensity (1.1) in a geometrically
clear way:
   it turns out that
on a $(n-1.n-1)$-dimensional supersurface embedded
in a $(n.n)$-dimensional odd symplectic superspace
 endowed with a volume form
 one can define in a natural way the
 odd semidensity
 whose values are
 odd vectors, symplectoorthogonal to this
  supersurface,
      and the
"truncated divergence" of this vector-valued semidensity
    is invariant semidensity
    (1.1). It has to be noted that our formula
  for this
   semidensity is very similar
   to the formula for the density corresponding
  to the mean curvature of the $1$-codimensional
   surface in the Euclidean space.
             $$ $$
      \centerline {\bf Section 2. Odd symplectic superspace.}
\medskip

         Let $E^{n.n}$  be a superspace with coordinates
        ($z^A=$  $x^1,\dots,x^n$, $\theta^1,\dots,\theta^n$).
        We say that this superspace is odd symplectic superspace
       if it is endowed with odd symplectic structure, i.e., if an odd
			 closed non-degenerated $2$-form $$
		    \Omega=\Omega_{AB}(z)dz^A dz^B\qquad
                    (p(\Omega)=1,\quad d\Omega=0)
                                                \eqno(2.1)
                          $$
                is defined on it [2,3] ($p$ is a parity:
         $p(x^i)=0,p(\theta^j)=1$).
        To the differential form (2.1) on the superspace
      $E^{n.n}$
      one can relate a function which
     for every point\footnote{$^*$}{More precisely, a point
  of superspace $E^{n.n}$ is $\Lambda$-point---
 $2n$--plet $(a^1,\dots,a^n,\alpha^1,\dots,\alpha^n)$ where
 $(a^1,\dots,a^n)$ are arbitrary even and $(\alpha^1,\dots,\alpha^n)$ are
arbitrary odd elements of an arbitrary Grassman algebra $\Lambda$.
  (We use the most general definition
 of superspace suggested by A.S. Schwarz
 as the functor on the category of Grassman
 algebras.[12])}
       defines the following skewsymmetric (in a supersense)
     odd bilinear
     form on  tangent vectors:
                         $$
                       \eqalign
                           {
            p\left(\Omega ({\bf X},{\bf Y})\right)&=
                   1+p({\bf X})+p({\bf Y})\,,
                           \cr
                   \Omega ({\bf X},{\bf Y})&=
                 -\Omega ({\bf Y},{\bf X})
                 (-1)^{p({\bf X})p(\bf Y)},
                         \cr
               \Omega ({\bf X},\lambda{\bf Y})&=
               \Omega ({\bf X}\lambda,{\bf Y}),
                        \cr
                       \Omega
         (\lambda{\bf Y}+\mu{\bf Z},{\bf X})&=
          \lambda\Omega({\bf Y},{\bf X})
                   +
          \mu \Omega({\bf Z},{\bf X})\,.
                         }
                                       \eqno (2.2)
                        $$
           In the coordinates:
                        $$
                     \Omega_{AB}=
          -\Omega_{BA}(-1)^{p(A)p(B)}=
                  \Omega\left(
            {\partial\over\partial z^A},
             {\partial\over\partial z^B}\right)\,,
        \left( p(\Omega_{AB})=1+p(A)+p(B)\right)\,.
                                                  \eqno (2.3)
                        $$

           Here ${\bf X}$, ${\bf Y}$, ${\bf Z}$
       are the vector fields
           $ X^A(z){\partial\over\partial z^A}$,
         $ Y^A(z){\partial\over\partial z^A}$,
       $Z^A(z){\partial\over\partial z^A}$
          (the left derivations of  functions on $E^{n.n}$).

(Differential form is usually considered as an element of an algebra
generated by  $z^A$ and $dz^A$, where the parity of $dz^A$ is
opposite to the parity of $z^A$. In this case
 $\Omega_{AB}$ $=\Omega_{BA}(-1)^{(p(A)+1)(p(B)+1)}$ instead of (2.3).
 The slight difference is eliminated by the transformation
 $\Omega_{AB}$ $\rightarrow$  $\Omega_{AB}(-1)^{B}$).

               From (2.2, 2.3) it follows that
                         $$
                      \Omega \left(
             X^A{\partial\over\partial z^A},
             Y^B{\partial\over\partial z^B}
                      \right)
                         =
                  X^A\Omega_{AB}Y^B
                       (-1)^
              { p(Y)p(B)+p(Y) }.
                                          \eqno (2.4)
                         $$

      In the same way as in the standard symplectic calculus one can
 relate to the odd symplectic structure (2.1) the odd Poisson
      bracket (Buttin bracket) [3]
                       $$
                     \{f,g\}=
          {\partial f\over\partial z^A}
                     (-1)^{p(f)p(A)+p(A)}
                          \Omega^{AB}
              {\partial g\over\partial z^B}
                                             \eqno (2.5)
                       $$
     where
                      $$
                 \Omega^{AB}
            = -\Omega^{BA}(-1)^{(p(A)+1)(p(B)+1)}=
               \{z^A,z^B\}
                     $$
 is the inverse matrix to $\Omega_{AB}$ :
                     $$
           \Omega^{AC}{ \Omega_{CB}}=
                   \delta^A_B\,.
                          \eqno (2.6)
                      $$
 To a function $f$ via (2.5) there corresponds
  the Hamiltonian vector field
                       $$
         {\bf D}_f=\{f,z^A\}{\partial\over\partial z^A}
             \quad{\rm and}\quad
         {\bf D}_f g=\{f,g\},\,
          \Omega({\bf D}_f,{\bf D}_g)=\{f,g\}\,.
                                     \eqno (2.7)
                      $$
    The condition of the closedness of the form (2.1) leads to the Jacoby
    identities:
                       $$
                      \eqalign
                       {
           \{f,\{g,h\} \}(-1)&^{(p(f)+1)(p(h)+1)}+
            \{g,\{h,f\} \}(-1)^{(p(g)+1)(p(f)+1)}+
                         \cr
            &\{h,\{f,g\} \}(-1)^{(p(h)+1)(p(g)+1)}=0\,.
                        }
                                       \eqno (2.8)
                       $$
    Using the analog of Darboux Theorem [13]
   one can consider the coordinates in
    which the symplectic structure (2.1) and the
   corresponding Buttin bracket
    have the canonical expressions.
    We call the coordinates
   $w^A$=$(x^1,\dots,x^n,\theta^1,...,\theta^n)$
  Darboux coordinates
    if in these coordinates holds
                      $$
            \Omega=I_{AB}dw^A dw^B \colon\quad
            \Omega\left(
     {\partial\over\partial x^i},{\partial\over\partial x^j}
            \right)=0\,,
           \Omega\left(
     {\partial\over\partial \theta^i},
 {\partial\over\partial \theta^j}\right)=0\,,
      \Omega\left({\partial\over\partial x^i},
      {\partial\over\partial \theta^j}\right)=-\delta^i_j\,,
                                         \eqno (2.9)
                      $$
  respectively
                        $$
                      \{f,g\}=
                      \sum_{i=1}^n
                       \left(
                 {\partial f\over\partial x^i}
                 {\partial g\over\partial \theta^i}
                           +(-1)^{p(f)}
                {\partial f\over\partial \theta^i}
                {\partial g\over\partial x^i}
                       \right)\,.
                                           \eqno (2.10)
                         $$

    Now on the odd symplectic superspace $E^{n.n}$ endowed
 with an odd symplectic structure
 (2.1) we consider  the
    following geometrical constructions:
     Let  $\bf\Psi$ be an odd nondegenerated vector field
 symplectoorthogonal to itself:
                  $$
             \Omega(\bf\Psi,\bf\Psi)=0
                                             \eqno (2.11)
                   $$
where
                      $$
                 {\bf \Psi}=
                   \sum_{i=1}^n
         \left( \Psi^i{\partial\over\partial\theta^i}
    				+
    	\Phi^i{\partial\over\partial x^i}\right)
                \quad \left(
            p(\Psi^i)=0,p(\Phi^i)=1
                    \right)
                                                    \eqno (2.12)
                       $$
             and at least for some index $i_0$
         the coefficient $\Psi^{i_0}$ is not nilpotent.
         (It is easy to see that this condition is invariant under
         the coordinate transformations.).

   (For example, to the even function $f$ such that
   $\{f,f\}=0$ there corresponds the Hamiltonian vector field
   ${\bf D}_f$ defined by (2.7), subject to
   condition (2.11) and to condition (2.12)
    if the gradient of $f$ is not nilpotent.)

   Let ${\bf\Pi}(z)$ be a field of  $(1.1)$-dimensional
   subspaces (planes)
  (${\bf\Pi} (z)\in TE_z^{n.n}$)
  which contain the vector field ${\bf\Psi} (z)$ and
 the symplectic structure induced on these planes is not
 degenerate.
  It means that
 there exists an even vector field
  ${\bf H}(z)$ such that
                       $$
    {\bf \Psi}(z)\,, {\bf H} (z)\in {\bf\Pi} (z)\quad
         {\rm and}\quad
         \Omega({\bf H} (z),{\bf\Psi} (z))=1.
                                                   \eqno (2.13)
                        $$
       To this field of planes ${\bf\Pi} (z)$ there corresponds the
       symplectoorthogonal projector $\hat \Pi (z)$ of the
       vectors in the tangent space to these planes:
                        $$
        {\hat \Pi}:\quad T_z E\rightarrow {\bf\Pi}(z),\quad
           {\hat\Pi}\Big\vert_{\bf\Pi}(z)={\bf id},\quad
                {\hat\Pi}{\bf X}=0\quad{\rm if}\quad
                       \Omega({\bf X},{\bf\Pi})=0\,.
                                          \eqno (2.14)
                         $$
      In the coordinates $\{z^A\}$ to the projector ${\hat \Pi}$
there corresponds the matrix-valued function $\Pi_A^B (z)$:
                         $$
          \hat\Pi\left({\partial\over\partial z^A}\right)=
           \Pi_A^B{\partial\over\partial z^B},\quad
                         {\rm so}\quad
                     {\hat \Pi}(X^A{\partial\over\partial z^A})=
                     X^B\Pi^A_B{\partial\over\partial z^A}
                     \quad{\rm and}\quad
                     \Omega({\bf X}, {\hat\Pi \bf Y})=
                     \Omega({\bf\Pi X, Y})\,.
                                            \eqno (2.15)
                          $$
  Later on we call
 $({\bf\Pi}(z),{\bf\Psi}(z))$
 or equivalently
 $({\hat\Pi}(z),{\bf\Psi}(z))$  an {\it odd normal pair}
if  ${\bf\Psi}(z)$ and
 ${\bf\Pi}(z)$  are defined by (2.11---2.15).
                 $$ $$
      \centerline {\bf Section 3. Special geometrical construction.}
\medskip
 In this section for odd normal pair
 $({\bf\Pi}(z),{\bf\Psi}(z))$ in an odd symplectic superspace
we construct a first-order divergence-like
differential operator which
 transforms it to a function on this superspace.

 Let in a superspace $E$, ${\bf X}(z)$ and ${\hat L}(z)$
 be a vector field and a linear operators field
 defined on $T_z E$ respectively. If ${\{z^A\}}$ are arbitrary
coordinates then for the pair
${({\hat L},{\bf X})}$
one can consider the function which depends on the
coordinate system $\{z^A\}$:
                          $$
            \partial({\hat L},{\bf X})^{\{z\}}=
              {\partial X^A(z)\over\partial z^B}
                     L^B_A(z)
                        (-1)^{p(X)p(B)+p(B)}\,.
                                                   \eqno (3.1)
                          $$
  Expression (3.1) is invariant under linear transformations of the
 coordinates $\{z^A\}$. In the general case
if ${\{w^A\}}$ and $\{z^A\}$ are two different coordinate systems
on $E$ then for the pair ${({L,\bf X})}$ we consider
                          $$
         \Gamma({\hat L},{\bf X})^{\{w\}}_{\{z\}}=
         \partial({\hat L},{\bf X})^{\{w\}}-
          \partial({\hat L},{\bf X})^{\{z\}}\,.
                                      \eqno(3.2)
                          $$
 From (3.1) and (3.2) it follows that
                           $$
        \Gamma({\hat L},{\bf X})^{\{w\}}_{\{z\}}=
                      X^Q(z)
                  \Gamma^A_{QB}
                (z\,\big\vert \{w\},\{z\})
                       L^B_A(z)
                        (-1)^{p(B)}
                           $$
where
                            $$
                       \Gamma^A_{BC}
            (z\,\big\vert{\{w\}},{\{z\}})=
        {\partial^2 w^K(z)\over\partial z^B \partial z^C }
        {\partial z^A(w)\over\partial w^K}\,.
                                         \eqno (3.3)
                            $$
 (In (3.3) the components of ${\hat L}$ and ${\bf X}$
are in the coordinates $\{z^A\}$.)

 \noindent  From definition
  (3.2) of the
  $ \Gamma(L,{\bf X})^{\{w\}}_{\{z\}}$ it follows that
  for
 three different coordinate systems
$\{w^A\}$, $\{z^A\}$ and $\{u^A\}$
                      $$
       \Gamma({\hat L},{\bf X})^{\{w\}}_{\{z\}}+
       \Gamma({\hat L},{\bf X})^{\{z\}}_{\{u\}}+
     \Gamma({\hat L},{\bf X})^{\{u\}}_{\{w\}}=0\,.
                                           \eqno (3.4)
                       $$
 Let ${\cal F}$ be a class of some coordinate systems
such that for a given pair $({\hat L},{\bf X})$
                           $$
  \forall \{z\},\{w\}\in {\cal F}\quad
          \Gamma({\hat L},{\bf X})^{\{w\}}_{\{z\}}=0\,.
                                           \eqno (3.5)
                           $$
Then to the class ${\cal F}$ one can relate
the first-order differential operator ${\cal D}$:
                          $$
  {\cal D}({\hat L},{\bf X})=
       {\bf\partial}({\hat L},{\bf X})^{\{z\}}+
[6~          \Gamma({\hat L},{\bf X})^{\{w\}}_{\{z\}}
                                       \eqno (3.6)
                        $$
where  $\{z\}$ are arbitrary coordinates on $E$ and $\{w\}$
are arbitrary coordinates from the
 class  ${\cal F}$.
 From (3.2, 3.5) it follows that the r.h.s. of (3.6)
 does not
 depend on the choice of these coordinates.

Before going to the considerations for an odd symplectic superspace
we consider an example where we come to the standard definition
of the divergence in superspace using (3.1---3.6).

 {\bf Example 1}. Let $E$ be a superspace with a volume form
 $d{\bf v}$ which
in coordinates
$\{w_0^A\}=\{x^1,\dots,x^n,\theta^1,\dots,\theta^m\}$ on $E$
is equal to
                           $$
d{\bf v}=dx^1\dots dx^n d\theta^1\dots d\theta^m.
                                   \eqno (3.7)
                           $$

We define ${\cal F}$ as a class of coordinate systems
in which the volume form $d{\bf v}$ is given by (3.7):
                          $$
    {\cal F}=\{\,\{w\}\colon\,Ber\left(
        {\partial w\over\partial w_0}
                \right)=1\,\}\,.
                                  \eqno (3.8)
                         $$
 ($Ber A$ is the superdeterminant of $A$.)
 It is easy to see that in this case
for arbitrary coordinates $\{z^A\}$
if ${\hat L}={\bf id}$ is identity operator and ${\bf X}$ is
an arbitrary vector field then
                     $$
                   \eqalign
                       {
         &\Gamma({\hat L},{\bf X})^{\{w\}}_{\{z\}}=
                      X^Q(z)
   {\partial^2 w^K(z)\over\partial z^Q \partial z^A }
        {\partial z^A(w)\over\partial w^K}
                  (-1)^{p(A)}=
                      \cr
           &X^Q(z){\partial\over\partial z^Q}
                      \log\left(
             Ber \left({\partial w\over\partial z}\right)
                      \right)=
                       X^Q(z)
       {\partial\log\rho(z)\over\partial z}
                           }
                                         \eqno (3.9)
                            $$
        where $\{w^A\}$ are arbitrary coordinates from the class
(3.8) and $\rho(z)dz^1\dots dz^{m+n}$ is the volume form (3.7)
in the coordinates $\{z^A\}$.
 The condition (3.5) is fulfilled
  and we come to the standard definition
 of the divergence in
a superspace endowed with volume form.
 For the pair
 $({\hat L},{\bf X})$ where ${\hat L}={\bf id}$
 and for class (3.8) the operator ${\cal D}({\hat L},{\bf X})$
 is the divergence of the
vector field ${\bf X}$ corresponding to the volume form $d{\bf v}$:
                     $$
                     \eqalign
                         {
     &{\cal D}({\hat L},{\bf X})=
       {\bf\partial}({\hat L},{\bf X})^{\{z\}}+
      \Gamma({\hat L},{\bf X})^{\{w\}}_{\{z\}}
                        =\cr
       &{\partial X^A(z)\over\partial z^A}
            (-1)^{p(X)p(A)+p(A)}+
        X^A {\partial\log\rho(z)\over\partial z^A}=
            div_{d{\bf v}}{\bf X}.
                        }
                                                \eqno (3.10)
                  $$
\medskip
 Now we return to the considerations of Section 2.

For the superspace $E^{n.n}$ endowed with the odd symplectic structure
 we consider a field
 $({\hat\Pi}(z),{\bf\Psi}(z))$,
 where  $({\hat\Pi}(z),{\bf\Psi}(z))$ is
 an odd normal pair in a vicinity of some
 point.
 (See
the end of the previous Section).
We denote by ${\cal F}_D$ the class of  Darboux coordinates
(2.9, 2.10) on $E^{n.n}$
 and apply constructions
 (3.1---3.6) in this case.

   {\bf Lemma.} If $({\hat\Pi},{\bf\Psi})$ is an odd normal
  pair in $E^{n.n}$ then
                      $$
\forall \{w\}\,,\{{\tilde w}\}\in {\cal F}_D,
\quad
  \Gamma({\hat\Pi},{\bf\Psi})^{\{{\tilde w}\}}_{\{w\}}=0\,.
                                  \eqno (3.11)
                       $$
 Using the statement of the Lemma we consider the
 action of the operator ${\cal D}_{can}$
 corresponding to the class
${\cal F}_D$ of Darboux coordinates by (3.6),
on the odd normal pair   $({\hat\Pi}(z),{\bf\Psi}(z))$:
                    $$
             {\cal D}_{can}
         ({\hat\Pi},{\bf\Psi})=
                   \left(
             {\partial\Psi^A(z)\over\partial z^B}
                      +
                  \Psi^Q(z)
       {\partial^2 w^K(z)\over\partial z^Q \partial z^B }
        {\partial z^A(w)\over\partial w^K}
                       (-1)^{p(B)}
                       \right)
                       \Pi^B_A(z)
                                   \eqno (3.12)
                          $$
where $\{w\}$ are arbitrary Darboux coordinates.
 From (3.6) and the Lemma follows

{\bf Theorem.} For the odd normal pair
  $({\hat\Pi}(z),{\bf\Psi}(z))$,
 ${\cal D}_{can}({\hat\Pi},{\bf\Psi})$ is an invariant
  geometrical object.

\noindent In particular if $\{w\}$ are
  arbitrary Darboux coordinates then
                         $$
             {\cal D}_{can}({\hat\Pi},{\bf\Psi})=
             {\partial\Psi^A(w)\over\partial w^B}
                \Pi^B_A(w)
                                              \eqno (3.13)
                          $$
does not depend on the choice of the Darboux coordinates $\{w\}$.

 Before proving the  Lemma we will consider

{\bf Example 2}.
 Let $E^{1.1}$ be a $(1.1)$-dimensional superspace endowed
with odd symplectic structure (2.1). Let
$w=(x,\theta)$ be some Darboux coordinates on it:
                         $$
            \{x,\theta\}=1,\,\{x,x\}=0\,.
                                       \eqno (3.14)
                         $$
It is easy to see that in this case the odd vector field
obeying to (2.11, 2.12)
is of the form
                        $$
               {\bf\Psi}=\Psi(x,\theta){\partial\over\partial\theta}
                                         \eqno (3.15)
                         $$
where $\Psi(x,\theta)$ is non-nilpotent even function.
The projector operator (2.14) is evidently
identity operator. So a normal pair is of the form
 $({\bf id},\Psi(x,\theta)\partial_\theta)$ and
                      $$
    {\cal D}_{can}({\bf id},{\bf\Psi})
    ={\partial \Psi(x,\theta)\over\partial\theta}\,.
                                           \eqno (3.16)
                      $$

 It is easy to see
that if ${x^{\prime},\theta^{\prime}}$ are some other
Darboux coordinates
 then they are related with coordinates $x,\theta$ by
canonical transformation
                          $$
                          \eqalign
                           {
  x^\prime&= f(x) \cr
  \theta^\prime&={\theta\over {df(x)/dx}}+\beta(x)
                           }
                                          \eqno (3.17)
                          $$
where $f(x)$ and $\beta(x)$ are even and odd functions
on $E^{1.1}$ respectively. (To obtain (3.17) from (3.14) one has to
note that in $E^{1.1}$ $\{x^{\prime},x^{\prime}\}=0$ $\rightarrow$
$x^{\prime}_\theta=0$.) It is easy to see from (3.17) that $$
   {\bf\Psi}=\Psi(x,\theta){\partial\over\partial\theta} =
    {\Psi(x,\theta)\over{df(x)/dx}}
    {\partial\over\partial\theta^\prime}
                                         \eqno (3.18)
                         $$
 and (3.16) does not change under
transformation (3.17), so in this case the
 statements of the Lemma and of the Theorem hold.

 Indeed, for this case one can say more about
(3.16). Let $(y,\eta)$ be an arbitrary coordinates
 on the $E^{1.1}$ and
a volume form $d{\bf v}$ on $E^{1.1}$ is defined by the
equation
                $$
     d{\bf v}={dy d\eta\over\{y,\eta\}}\,.
                              \eqno (3.19)
               $$
        Then one can check using
(3.10) that ${\cal D}_{can}({\bf id},{\bf\Psi})$ in (3.16) is the
divergence of the vector field ${\bf\Psi}$ by the volume form
   (3.19) and does not depend on
the choice of  coordinates $(y,\eta)$.

 Now we prove the Lemma. Let
                          $$
      z^A=w^B L^A_B+w^Bw^CT^A_{BC}+o(w^2)
                                      \eqno (3.20)
                        $$
be arbitrary canonical transformations from Darboux coordinates
 $\{z\}$ to Darboux coordinates $\{w\}$
in a vicinity of the point $z=0$.
For proving the Lemma we have to show that for transformation
 (3.20)
                       $$
       \Gamma({\hat\Pi},{\bf \Psi})^{\{w\}}_{\{z\}}
                   \Big\vert_{z=0}=
                  \Psi^Q(z)
                  {\partial^2 w^K(z)\over\partial z^Q z^B }
             {\partial z^A(w)\over\partial w^K}
                   \Pi^B_A(-1)^B\Big\vert_{z=0}
                        =0.
                                          \eqno(3.21)
                       $$
(We can consider (3.20) without loss of generality, since
in Darboux coordinates the translation is
 obviously the canonical transformation.)
We include the Darboux transformation (3.20) in the chain of
Darboux transformations:
                   $$
   \{z\}{\buildrel {\rm linear}\over\longrightarrow}
   \{{\tilde z}\}\longrightarrow \{w\}
   {\buildrel{\rm linear}\over\longrightarrow }\{{\tilde w}\}
                                                  \eqno (3.22)
                     $$
which obey to the following conditions:

  a) The transformation
                       $$
              \{z\}\rightarrow \{{\tilde z}\}
                                            \eqno (3.23)
                        $$
is the linear canonical transformation
such that in the  Darboux coordinates ${\tilde z}$
 $=({\tilde x}^i,{\tilde\theta}^k)$
                         $$
                      {\hat\Pi}
                 {\partial\over\partial {\tilde x^1}}=
                     {\partial\over\partial {\tilde x^1}},
[6~                    \quad
              {\hat\Pi}
        {\partial\over\partial {\tilde \theta^1}}=
               {\partial\over\partial {\tilde \theta^1}}\,.
                                                 \eqno (3.24)
				 $$
b) The transformation
                          $$
               \{w\}\longrightarrow \{{\tilde w}\}
                                                   \eqno (3.25)
                        $$
is the linear canonical transformation such that
                       $$
            {\tilde w}^A= {\tilde z}^A+ o({\tilde z})\,.
                                                 \eqno (3.26)
                     $$
   From (3.2, 3.4) it follows that
                    $$
     \Gamma({\hat\Pi},{\bf \Psi})^{\{w\}}_{\{z\}}=
\Gamma({\hat\Pi},{\bf \Psi})^{\{w\}}_{\{{\tilde w}\}}+
\Gamma({\hat\Pi},{\bf \Psi})^
{\{{\tilde w}\}}_{\{{\tilde z}\}}+
\Gamma({\hat\Pi},{\bf \Psi})^{\{{\tilde z}\}}_{\{z\}}\,.
                                        \eqno (3.27)
                    $$
 But
 $\Gamma({\hat\Pi},{\bf \Psi})^{\{w\}}_{\{{\tilde w}\}}$
 and
 $\Gamma({\hat\Pi},{\bf \Psi})^{\{{\tilde z}\}}_{\{z\}}$
 in (3.27) are zero because the corresponding transformations
are linear.
  The transformation
  $\{{\tilde w}\}\longrightarrow\{{\tilde z}\}$ is not linear
 but from (3.24, 3.26, 2.11, 2.12) it follows that
 ${\bf\Psi}^{\{{\tilde z}\}}\big\vert_{z=0}=
  \Psi\partial_{{\tilde\theta}^1}$
 and
                            $$
           \Gamma({\hat\Pi},{\bf \Psi})^
          {\{{\tilde w}\}}_{\{{\tilde z}\}}
             \Big\vert_{z=0}=
                  \Psi
                   {
        \partial^2 x^{\prime 1}
              \over
  \partial {\tilde\theta}^1\partial {\tilde x}^1
                   }
                   = 0
                                           \eqno (3.28)
                    $$
 because $\{x^{\prime 1},x^{\prime 1}\}=0$.
($\{{\tilde w}\}=(x^{\prime 1},\cdots)$).
(Compare with Example 2).
 So (3.21) is obeyed.

\noindent Lemma is proved.

    Considering a pair $({\hat L}(z),{\bf X}(z))$
 we constructed in this section the divergence-like operator
(3.12, 3.13) in an odd symplectic superspace in the case when
${\hat L}$ is a projector on a $(1.1)$-dimensional symplectic
subspace and ${\bf X}$ is an odd vector
in it which is symplectoorthogonal to itself. In the case of even
symplectic structure this construction can be carried out
in a more general case and it is trivial
 because in this
case there exists a volume form corresponding to the
symplectic structure. Indeed if in a superspace
 $E^{2m.n}$ endowed with
even symplectic structure there is given
  a pair $({\hat L}(z),{\bf X}(z))$  where
  ${\hat L}(z)$ is a symplectoorthogonal projector on
   $(2p.q)$-dimensional symplectic subspaces in $T_z E^{2m.n}$
  and ${\bf X}(z)$ is an arbitrary vector field then
 to the class ${\cal F}$ of the Darboux
  coordinates on this superspace there
 corresponds  ${\cal D}({\hat L,{\bf X}})$ defined by (3.6) which
is a straightforward generalization of (1.4).

       \bigskip
      \centerline {\bf Section 4. Truncated divergence.}
\medskip

 We consider in this section the odd vector field
 which is defined on the points of $(n-1.n-1)$-dimensional
 nondegenerate supersurface embedded in the
 $(n.n)$-dimensional odd symplectic superspace
 with a volume form. In case of this vector field being
 symplectoorthogonal to this supersurface and to itself
  using the geometrical object
 ${\cal D}_{can}({\hat\Pi},{\bf\Psi})$ for an odd normal pair
 we define the first order differential operator
 on it (truncated divergence) whose action on this field
 gives the function on this supersurface.

 Let $M^{n-1.n-1}$ be an arbitrary nondegenerate
 $(n-1.n-1)$-dimensional supersurface
 embedded in a superspace
$E^{n.n}$ which is endowed with odd
symplectic structure (2.1 ) and the volume form $d{\bf v}$
                        $$
         d{\bf v}=\rho(z)dz^1\dots dz^{2n}\,.
                          \eqno (4.1)
                        $$
  (The supersurface $M$ is nondegenerate if the symplectic structure
  of $E^{n.n}$ induces nondegenerate symplectic structure on $M$.)
Let $z^A=z^A(\zeta^\alpha)$ be a local parametrization
of the supersurface $M$. (For us it will be convenient to denote by the
letter "$\alpha$" the coordinates on $M$).)
 In the coordinates the induced two-form
  $\Omega_{\alpha\beta}d\zeta^\alpha d\zeta^\beta$
 on $M$ is given by the following equation (See Section 2):
                    $$
     \Omega_{\alpha\beta}(z(\zeta))=
        \Omega\left(
                { \partial_\alpha z^A}
         {\partial\over\partial z^A},
         {\partial_\beta z^B}
         {\partial\over\partial z^B}\right)
                   =
             {\partial_\alpha z^A}
             \Omega_{AB}
              \partial_\beta z^B
               (-1)^
              {s(B,\beta)}         \,.
                                 \eqno (4.2)
                    $$
 Hereafter  we use notations
   $\partial_\alpha z^A={\partial z^A\over\partial \zeta^\alpha},
   \partial_\alpha f={\partial f\over\partial \zeta^\alpha},\dots$
    for derivatives along surface and
                              $$
      (-1)^{s(B,\beta)}=(-1)^{p(B))p(\beta)+p(\beta)}
                                              \eqno (4.2a)
                            $$
 for sign factor.

The induced Poisson bracket structure on $M$
  $\{\,,\,\}_M$ according to (2.5, 2.6)
  is defined by the matrix $\Omega^{\alpha\beta}$ which
 is inverse to the matrix ${\Omega}_{\alpha\beta}$.
 Using this induced symplectic structure one can construct the
symplectoorthogonal projector on $TM$

\noindent $T_{z(\zeta)}E^{n.n}
{\buildrel {\hat P}(z(\zeta))\over\longrightarrow}
T_{z(\zeta)}M^{n-1.n-1}$ which can be expressed in
terms of $\Omega_{AB}$ and $\Omega^{\alpha\beta}$
                     $$
                    {\hat P}:\,
                   P_A^B(z(\zeta))=
                   { \Omega}_{AK}(z(\zeta))\cdot
                \{z^K(\zeta),z^B(\zeta)\}_M=
               { \Omega}_{AK}(z(\zeta))
                  \partial_\alpha z^K(-1)^{s(K,\alpha)}
                \Omega^{\alpha\beta}
                 \partial_\beta z^B\,.
                                                      \eqno (4.3)
                        $$

 The operator
                           $$
              {\hat\Pi}={\bf id}-{\hat P}
                                \eqno (4.4)
                           $$
in every point $z(\zeta)$ is a symplectoorthogonal projector
on the $(1.1)$-dimensional subspace ${\bf\Pi}(z(\zeta))$
in $T_{z(\zeta)}E$
which is symplectoorthogonal  and transversal
to $T_{z(\zeta)}M$ (See 2.14, 2.15).
It is easy to see that there exists the odd vector field
${\bf\Psi}(z(\zeta))$ belonging to ${\bf\Pi}(z(\zeta))$
 which is defined on the points of
 the supersurface $M$, is symplectoorthogonal to $M$ and which is
symplectoorthogonal to itself and  non-degenerate
 (2.11, 2.12):
                        $$
                   {\hat P}{\Psi}=0,\,
          \Omega ({\bf\Psi},{\bf\Psi})=0.\,
             ({\rm In\,the\,components}\quad
           {\Psi}^A P_A^B= 0,\quad
           \partial_\alpha z^A P_A^B=\partial_\alpha z^B)\,.
                                      \eqno (4.5)
                        $$
(This vector field is fixed uniquely up to the multiplication
by  an even non-nilpotent function of $\zeta$.
 \noindent ${\bf\Psi}$ $\rightarrow$
 $f(z(\zeta)){\bf\Psi}$.)

 The field ${\bf\Psi}$ and the projector ${\hat \Pi}$
 form the odd normal  pair $({\hat\Pi},{\bf\Psi})$
in the points $z(\zeta)$ of the supersurface $M$.

   For example if the supersurface $M^{n-1.n-1}$ is defined
by equations
                  $$
            f(z)=0,\,\varphi(z)=0
                                       \eqno (4.6)
                  $$
  then for arbitrary point $z_0$ on this supersurface, the
 vectors  ${\bf D}_\varphi$ and
   ${\bf D}_f$ (see 2.7) are the basis vectors of
   the subspace ${\bf\Pi}(z_0)$ and the vector
                         $$
              {\bf \Psi}=
                   \left(
                 {\bf D}_f-
          {\{f,f\}\over 2\{f,\varphi\}}
                   {\bf D}_\varphi
                       \right)
                    \Big\vert_{z_0}
                                       \eqno (4.7)
                         $$
 is subject to the conditions  (4.5).

 Now for this odd vector field which is defined
only on the supersurface $M$ we will construct
the "truncated divergence". For this purpose
we consider an odd normal pair
 $({\tilde{\hat\Pi}},{\tilde{\bf\Psi}})$ in $E^{n.n}$
 in the vicinity of the arbitrary point $z(\zeta)$
 which is a prolongation of the odd normal pair
$({\hat\Pi},{\bf\Psi})$:
                   $$
         {\tilde{\bf \Psi}}\big\vert_{z=z(\zeta)}=
             {\bf \Psi}(z(\zeta)),\quad
          {\tilde {\hat\Pi}}\big\vert_{z=z(\zeta)}=
               {\hat\Pi}(z(\zeta))
                                    \eqno (4.8)
                       $$
   and define the truncated divergence in the following way:
                       $$
            Div_{trunc}{\bf\Psi}(z(\zeta))=
                    \left(
              div_{d{\bf v}}{\tilde{\bf\Psi}}-
                      {\cal D}_{can}
     ({\tilde {\hat\Pi}},{\tilde {\bf\Psi}})
                           \right)
                      \Big\vert_{z=z(\zeta)}\,.
                                   \eqno (4.9)
                       $$
 In (4.9) ${\cal D}_{can}({\tilde {\hat\Pi}},{\tilde {\bf\Psi}})$
is given by (3.12, 3.13) for the odd normal pair
 $({\tilde {\hat\Pi}},{\tilde {\bf\Psi}})$
and $div_{d{\bf v}}{\bf\Psi}$ is  divergence (3.10)
of the vector field corresponding to the volume
form $d{\bf v}$ defined by (4.1).

 One can see that the r.h.s. of equation
(4.9) indeed does not depend on the prolongation
 $({\tilde {\hat\Pi}},{\tilde {\bf\Psi}})$
of the odd normal
pair $({\hat\Pi},{\bf\Psi})$.

Using equations (4.3 -- 4.5) for the projectors
${\hat P}$ and ${\hat \Pi}$, formulae (3.9, 3.10) for
the divergence, formula (3.12) for
 ${\cal D}_{can}({\tilde {\hat\Pi}},{\tilde {\bf\Psi}})$
  and the fact that
                  $$
               \partial_\beta  z^B
     {\partial {\tilde \Psi}(z)^A\over\partial z^B}
             \Big\vert_{z=z(\zeta)}=
           \partial_\beta \Psi^A(z(\zeta))
                            \eqno (4.10)
                 $$
  depends only on ${\bf\Psi}(z(\zeta))$,
 we rewrite $Div_{trunc}{\bf\Psi}(z(\zeta)$
 in arbitrary coordinates $\{z^A\}$ in the following way
                    $$
		   \eqalign
		     {
               Div_{trunc}{\bf\Psi}(z(\zeta))&=
                      \cr
                     \left(
               \partial_\beta  {\Psi}^A
                       +
                \partial_\beta  z^B\Psi^Q
       {\partial^2 w^K(z)\over\partial z^Q \partial z^B }
        {\partial z^A(w)\over\partial w^K}
                       (-1)^{p(B)}
                       \right)&
                     \Omega_{AK}
                    \partial_\alpha z^K
                 \Omega^{\alpha\beta}
                (-1)^{s(K,\alpha)}
			\cr
                         +
         \Psi^A{\partial \over\partial z^A}
                       \left(
                     \log \rho(z)-
                         \log
                          \left(
            Ber{\partial w\over\partial z}
                           \right)
                   \right)& \Big\vert_{z=z(\zeta)}
                          }
                                       \eqno (4.11)
                        $$
  where $\{w^A\}$ are any  Darboux coordinates.

 In arbitrary Darboux coordinates $\{w^A\}$
using formula (3.13) for
 ${\cal D}_{can}({\tilde {\hat\Pi}},{\tilde {\bf\Psi}})$
 we arrive at the following expression for
     $Div_{trunc}{\bf\Psi}(w(\zeta))$:
                    $$
     Div_{trunc}{\bf\Psi}(w(\zeta))=
               \partial_\beta\Psi^A
                I_{AK}
             \partial_\alpha w^K
       \Omega^{\alpha\beta}(w(\zeta))
       (-1)^{s(K,\alpha)}
                  +
               \Psi^A
    {\partial \log\rho(w)\over\partial w^A}
         \Big\vert_{w=w(\zeta)}\,.
                                    \eqno (4.12)
                   $$

Indeed  $Div_{trunc}{\bf\Psi}$ depends
 only on the values of ${\bf\Psi}$ and it
does not depend on the derivatives $\partial_\alpha\Psi^A$
because ${\bf\Psi}$ is symplectoorthogonal to $M^{n-1.n-1}$

\noindent
($\Psi^A\Omega_{AK}\partial_\alpha z^K(-1)^{s(K,\alpha)}=0$).
  One can rewrite for example (4.12) in the following way:
                    $$
                \eqalign
                    {
     Div_{trunc}{\bf\Psi}(w(\zeta))=
                  \cr
                  \Psi^A
                  \left(
                -I_{AK}
            \partial_\beta\partial_\alpha w^K
       \Omega^{\alpha\beta}(w(\zeta))
       (-1)^{s(K,\alpha+\beta)+p(\beta)}
                  +
    {\partial \log\rho(w)\over\partial w^A}
         \Big\vert_{w=w(\zeta)}
                    \right)\,.
                      }
                                    \eqno (4.12a)
                     $$
(For the definition of $s(K,\alpha+\beta)$ see (4.2a).)

  In the case if there exist Darboux coordinates
 $\{w\}$ in which the volume form $d{\bf v}$ is trivial,
                $$
   d{\bf v}=dz^1\dots dz^{2n}\quad (\rho=1)
                                   \eqno (4.13)
               $$
then
                   $$
     Div_{trunc}{\bf\Psi}(w(\zeta))=
                  -\Psi^A
                 I_{AK}
             \partial_\beta\partial_\alpha w^K
             (-1)^{s(K,\alpha+\beta)+p(\beta)}
       \Omega^{\alpha\beta}(w(\zeta))\,.
                                    \eqno (4.14)
                   $$

 Equation (4.14) defines  $Div_{trunc}{\bf\Psi}(w(\zeta))$
in arbitrary Darboux coordinates (if they exist)
which are subject to condition (4.13).

 \bigskip
    \centerline {\bf Section 5.  The odd invariant semidensity.}
  \medskip
 Now we are well prepared for writing the formula for
 odd invariant semidensity using  constructions (4.9---4.12) for
 truncated divergence.
 The constructions of this section are founded
 on the following remark.
 Let $M^{n-1.n-1}$ be an arbitrary nondegenerate
 $(n-1.n-1)$-dimensional supersurface
 embedded in the odd symplectic superspace $E^{n.n}$ and a field
 ${\bf\Psi}$ on this supersurface obeys to conditions (4.5).
 Then the r.h.s. of (4.11, 4.12) by its definition is invariant
 under coordinate transformations of the superspace $E^{n.n}$
 and does not depend on the parametrization $z^A=z^A(\zeta^\alpha)$
 of the supersurface $M^{n-1.n-1}$.
 So the equations (4.11, 4.12) define a density of the weight
 $\sigma=0$ which is defined on
 $(n-1.n-1)$-dimensional supersurfaces.
  Moreover,
 one can see from (4.12a) that if ${\bf\Psi}$ is
 a density of an arbitrary
weight $\sigma$ which is defined on supersurfaces
  $M^{n-1.n-1}$  and takes values
in the odd vector fields obeying to (4.5) then the
truncated divergence of this density is
 the density of the
weight $\sigma$ which is defined on $(n-1.n-1)$-dimensional
nondegenerate supersurfaces and takes numerical values.

 Let, as in Section 4, $M^{n-1.n-1}$ be an
 arbitrary nondegenerate
 supersurface in an odd symplectic superspace
 $E^{n.n}$ which is endowed with
 symplectic structure (2.1) and
 volume form (4.1).
 Now we will construct the semidensity
 on the $M$ which takes values
 in the odd vectors ${\bf\Psi}$ obeying to conditions (4.5).
 Let the vectors
 $({\bf e}_1,...,{\bf e}_{n-1};{\bf f}_1,...,{\bf f}_{n-1})$
 constitute
 a basis of the tangent space $T_{z(\zeta)}$ in arbitrary point
 $z(\zeta)$ of the supersurface $M^{n-1.n-1}$
 ( ${\bf e}_i$ are even vectors and ${\bf f}_i$ are odd ones).
 Let ${\bf\Psi}(z(\zeta))$ and ${\bf H}(z(\zeta))$
 be respectively an odd and an even vector fields
 which belong to ${\bf\Pi}(z(\zeta))$ (see (4.5))
 such that $({\bf\Pi}(z(\zeta)),{\bf\Psi}(z(\zeta))$
form an odd normal pair
                    $$
          \Omega ({\bf \Psi},{\bf \Psi})=0\quad{\rm and}\quad
                \Omega ({\bf H},{\bf \Psi})=1\,.
                                     \eqno (5.1)
                   $$
These conditions fix the vector fields
 ${\bf H}$ and ${\bf\Psi}$ up to the
transformation
              $$
 {\bf H}\rightarrow {1\over\lambda} {\bf H}+\beta {\bf\Psi},
\quad {\bf\Psi}\rightarrow \lambda{\bf\Psi}
                                          \eqno (5.2)
              $$
where $\lambda$ is an arbitrary even function
(taking values in non-nilpotent numbers)
and $\beta$ is an arbitrary odd function
(compare with 3.17). Using (5.2)
one can choose the vector field ${\bf\Psi}$ (but not the vector
field {\bf H}) in the unique way by imposing the
normalization condition via volume form (4.1):
                   $$
               d{\bf v}
   ({\bf e}_1,...,{\bf e}_{n-1},{\bf H};
   {\bf f}_1,...,{\bf f}_{n-1},{\bf\Psi})=1
                             \eqno (5.3)
                   $$
We arrive at the function
               $$
       {\bf\Psi}=
      {\bf\Psi}(z(\zeta),
   {\bf e}_1,...,{\bf e}_{n-1};{\bf f}_1,...,{\bf f}_{n-1})\,.
                                \eqno (5.4)
              $$
 which depends on points
 $z(\zeta)$ of the supersurface $M^{n-1.n-1}$ and the bases
 $({\bf e}_1,...,{\bf e}_{n-1};$ ${\bf f}_1,...,{\bf f}_{n-1})$
 in the $T_{z(\zeta)}M^{n-1.n-1}$  and which takes values in the
 odd vector fields obeying to condition (4.5).
 This function is defined uniquely
 by  condition (5.1) and by
 normalization condition (5.3).
 It is easy to see that under the change of the basis
 the function (5.4) is multiplied by the
 {\it square root of the corresponding Berezinian}
 \footnote{$^*$}{It is interesting to note that
 these considerations for obtaining the formula
 for invariant vector--valued semidensity
 are similar to
 the considerations for obtaining the formula
 for the invariant density on the
 lagrangian surfaces in $E^{n.n}$ suggested by A.S.
 Schwarz [5].}.
 For example
 if ${\bf e_1}$ $\rightarrow$  $\lambda{\bf e_1}$
 and ${\bf f_1}$ $\rightarrow$  $\mu{\bf f_1}$
 then ${\bf \Psi}$ $\rightarrow$
 $\sqrt{{\lambda\over\mu}}{\bf \Psi}.$

 If $z^A=z^A(\zeta^\alpha)$ is any parametrization
 of the supersurface $M^{n-1.n-1}$ where
 $\{\zeta^\alpha\}$
 $=(\xi^1,\dots,\xi^{n-1};\nu^1,\dots,\nu^{n-1})$
 are even and odd  parameters of this supersurface
 then  considering as the basis vectors
                  $$
             {\bf e}_1
                  =
   {\partial z^A\over\partial\xi^1}
   {\partial\over\partial z^A},
          \dots,
      {\bf e}_{n-1}
            =
   {\partial z^A\over\partial\xi^{n-1}}
   {\partial\over\partial z^A};\,
             {\bf f}_1
                =
   {\partial z^A\over\partial\nu^1}
   {\partial\over\partial z^A},
          \dots,
             {\bf f}_{n-1}
                 =
   {\partial z^A\over\partial\nu^{n-1}}
    {\partial\over\partial z^A}
               $$
 we see that (5.4) defines odd vectors valued semidensity
 ${\bf\Psi}(z(\zeta),{\partial z\over\partial \zeta})$
of the rank $k=1$  on nondegenerate
 $(n-1.n-1)$-dimensional supersurfaces.
The truncated divergence of this
semidensity is the odd semidensity of the rank $k=2$.
Using the formula
 (4.12a) we arrive at this odd invariant semidensity:
                    $$
		 \eqalign
		     {
     &A\left(w(\zeta),{\partial w\over\partial\zeta},
  {\partial^2 w^A\over\partial\zeta\partial\zeta}\right)
                   =
           Div_{trunc}{\bf\Psi}
  \left(w(\zeta),{\partial w\over\partial \zeta}\right)=
                   \cr
                &\Psi^A
   \left(w(\zeta),{\partial w\over\partial \zeta}\right)
                \left(
                -I_{AK}
               \partial_\beta\partial_\alpha  w^K
         \Omega^{\alpha\beta}(w(\zeta))
           (-1)^{s(K,\alpha+\beta)+p(\beta)}
                  +
    {\partial \log\rho(w)\over\partial w^A}
         \Big\vert_{w=w(\zeta)}
                  \right)
		   }
                                    \eqno (5.5)
                   $$
 where $\{w^A\}$
 are arbitrary Darboux coordinates in the $E^{n.n}$
and $w(\zeta)$ is parametrization of supersurface $M^{n-1.n-1}$.
 In a case if there exist Darboux
coordinates $\{w^A\}$ in which
the volume form is trivial ($\rho=1$) the formula (5.5)
according to (4.14) is reduced to
                    $$
              A\left(w(\zeta),{\partial w\over\partial\zeta},
            {\partial^2 w\over\partial\zeta\partial\zeta}\right)=
                -\Psi^A
  \left(w(\zeta),{\partial w\over\partial \zeta}\right)
                I_{AK}
                  \partial_\beta\partial_\alpha w^K
          \Omega^{\alpha\beta}(w(\zeta))
           (-1)^{s(K,\alpha+\beta)+p(\beta)}\,.
                                    \eqno (5.6)
                   $$
 The semidensity (5.5) is nothing but the
 semidensity obtained in [8]. (See (1.1).)
 To compare (5.5) with  (1.1)
 we consider in the vicinity of arbitrary point $z_0$
 the Darboux coordinates
 in which the supersurface is flat up to the second order derivatives:
 $\{w^A\}=(x^1,\dots,x^n,\theta^1,\dots,\theta^n)$  and parameters
 $\{\zeta^\alpha\}=(\xi^1,\dots,\xi^{n-1},\nu^1,\dots,\nu^{n-1})$
 such that pa\-ra\-met\-ri\-zation
 $w^A=w^A(\zeta^\alpha)$ is of the form:
                       $$
                    \eqalign
                      {
              x^{n}&=o(\zeta^2),
                      \cr\
               \theta^n &=o(\zeta^2),
                      \cr
       x^i &=\xi^i+o(\zeta^2)\quad{\rm if}\quad 1\leq i\leq n-1,
                       \cr
       \theta^i &=\nu^i+o(\zeta^2)\,{\rm if}\quad 1\leq i\leq n-1\,.
                   }
                                           \eqno (5.7)
                         $$
 (In (5.7) $z_0=0$.)
 It is evident that the vector valued semidensity (5.4)
in the point $z_0$ in  parametrization (5.7)
 is equal to $\partial_{\theta^n}$ and
 $w^A_{\alpha\beta}\vert_{z_0}=0$, so semidensity (5.5) in the
point $z_0$ is equal to
              $$
       {\partial\log\rho\over\partial\theta^{n}}\Big\vert_{z=z_0}\,.
                                \eqno(5.8)
              $$
Returning to (1.1) we see that to (5.7) there correspond the functions
 $f=x^n+o(z^2)$ and $\varphi=\theta^n+o(z^2)$, so the dual density
  ${\tilde A}$ in (1.1) in this point is equal to
                 $$
          \Delta f\Big\vert_{z=z_0}=
            div_{d{\bf v}}{\bf D}_f\Big\vert_{z=z_0} =
       {\partial\log\rho\over\partial\theta^{n}}\Big\vert_{z=z_0}
               $$
and coincides (up to a constant) with (5.8).

\bigskip
    \centerline {\bf Section 6.  Discussions.}
  \medskip

The constructions of truncated divergence and of odd semidensity
have analogues in the standard differential geometry.

The constructions of Section 4 have evident analogues
 for the hypersurfaces in a Riemannian space.
For example to
 (4.14) there corresponds the following construction.
 Let $C$ be an $(n-1)$-dimensional surface (hypersurface)
 embedded in an $n$-dimensional Euclidean space $E^n$
 and ${\bf R}$ be a vector field defined on the surface
 $C$ which is orthogonal to $C$. Then one can consider
in analogy with (4.14) the function
                   $$
         R^i G_{ik}\partial_\beta\partial_\alpha x^k g^{\alpha\beta}
                                                     \eqno (6.1)
                 $$
which does not depend on the parametrization
$x^i=x^i(\xi^\alpha)$ of the surface $C$ and
on the choice of the Euclidean coordinates $\{x^i\}$ in $E^n$.
 In (6.1) $G_{ik}=\delta_{ik}$ is the metric tensor in $E^n$,
 $g_{\alpha\beta}=\partial_{\alpha}x^i G_{ik}
 \partial_{\beta}x^k_\beta$
 is the metrics induced on the surface $C$,
 $g^{\alpha\beta}$ $=(g)^{-1}_{\alpha\beta}$ is inverse metric tensor.

More interesting  is to compare
the odd semidensity constructed in Section 5 with the
mean curvature of the hypersurfaces in  Euclidean space.

In analogy with considerations of Section 5 one can consider
for the hypersurface $C$ in $E^n$  invariant
 density $R^i(x(\xi),{\partial x\over\partial\xi})$ of the rank $k=1$
 which takes
values in the vectors orthogonal to this surface.
  By these conditions it is fixed uniquely
(up to multiplication by the constant):
                  $$
      R^i(x(\xi),{\partial x\over\partial\xi})
                     =
         n^i(x(\xi))\sqrt{det (g_{\alpha\beta})}
                                               \eqno (6.2)
                    $$
where  $ n^i(x(\xi))$ is the unit vector field orthogonal to this
surface and  $\sqrt{det (g_{\alpha\beta})}$ is the density of
the volume form induced on the surface. Applying (6.1)
 to vector-valued density (6.2) we come
 in analogy with (5.6) to the following density
of the second rank
                  $$
    H(x(\xi),{\partial x\over\partial\xi},
         {\partial^2 x\over\partial\xi\partial\xi})=
         n^i(x(\xi))\partial_\alpha\partial_\beta x^i
       g^{\alpha\beta}\sqrt{det (g_{\alpha\beta})}\,.
                                              \eqno (6.3)
               $$
 It is easy to see that
 $n^i(x(\xi))\partial_\alpha\partial_\beta x^i$
are components of the second quadratic form for the hypersurface $C$
 and the density (6.3) corresponds to
the mean curvature [14].

 Here we want to note that in spite of the
fact that odd semidensity (5.5)
 can not be integrated over surfaces
 ($A^2=0$) one can consider  the equation
                    $$
                    A\equiv 0\,.
                                            \eqno (6.4)
                     $$
which extracts locally the class of
$(n-1.n-1)$ supersurfaces (in an odd
 symplectic superspace endowed with
a volume form) on which (6.4) is satisfied.

 For example, from (5.6---5.8) it follows that
to this class there belong supersurfaces which can be locally defined by
equations $x^n=\theta^n=0$
 if there exist Darboux coordinates

  \noindent
  $(x^1,\dots, x^n,\theta^1\dots,\theta^n)$
in which the volume form is trivial.

 The analogous condition for mean curvature (6.3)
                    $$
                    H\equiv 0\,
                                            \eqno (6.5)
                  $$
    is the solution of the
variational problem for minimal surfaces.
 Mean curvature (6.3) of the hypersurface $C$ in $E^n$
 is identically equal to zero iff the surface $C$
locally is extremal for the
 "surface " functional which is equal to the integral
over the surface of the volume element (density)
 $\sqrt {det g_{\alpha\beta}}$
corresponding to the metrics $g_{\alpha\beta}$
induced on the surface [14].
       \bigskip
       \centerline {\bf Acknowledgments}

 I want to express my deep gratitude to
A.S. Schwarz and D.A.Leites who got me acquinted
with the problems of odd symplectic geometry many years ago.

 I am deeply grateful to J.C. Hausmann,
 to Th. Voronov and to A.V.Karabegov
 who encouraged
 me for this work.
 \vfill
 \eject

\def \m {\medskip}
   \centerline {\bf References}
\m
[1]  Batalin,I.A., Vilkovisky,G.A.:  Gauge algebra and quantization.
 {\it Phys.Lett.} {\bf 102B}, (1981), 27--31.
\smallskip
   Closure of the gauge algebra,
 generalized Lie equations and Feynman rules.
  {\it Nucl.Phys.} {\bf B234}, (1984), 106--124.
\m
[2] Berezin,F.A.:
{\it Introduction to Algebra and Analysis with Anticommuting
 Va\-ri\-ables.}  Moscow, MGU (1983).
\medskip
 [3]. Leites,D.A.: The new Lie superalgebras and Mechanics.
   {\it Docl. Acad. Nauk SSSR} {\bf 236}, (1977), 804--807.
\smallskip
 Introduction to the theory of supermanifolds.
      {\it Usp. Mat. Nauk} {\bf 35}, No.1, 3--57.
\smallskip
  {\it The theory of Supermanifolds.}
   Karelskij Filial AN SSSR (1983).
\m
[4] Khudaverdian,O.M.:  Geometry of superspace provided by Poisson
brackets of different gradings. {\it  J. Math. Phys.}
 {\bf 32}, (1991), 1934--1937.
\m
[5]    Schwarz,A.S.: Geometry of Batalin--Vilkovisky Formalism.
  {\it Commun. Math. Phys.} {\bf 155}, (1993), 249--260.
\m
[6]  Khudaverdian,O.M., Nersessian,A.P.:  On the geometry of
Batalin--Vilkovisky formalism.
 {\it Mod. Phys. Lett.} {\bf A8}, (1993), 2377--2385.
 \smallskip
  Batalin--Vilkovisky Formalism and
  Integration Theory on Manifolds.
  {\it J. Math. Phys.} {\bf 37}, (1996), 3713--3721.
 \m
[7] Khudaverdian,O.M., Schwarz,A.S., Tyupkin,Yu.S.,
 Integral invariants for Supercanonical Transformations
 {\it Lett. Math. Phys.}, {\bf 5} (1981), 517--522.
\medskip
 [8] Khudaverdian,O.M., Mkrtchian,R.L.:
  Integral invariants of Buttin bracket.
  {\it Lett. Math. Phys.} {\bf 18}, (1989), 229--234.
\m
[9] Gayduk,A.V.,Khudaverdian,O.M, Schwarz,A.S,:
  Integration on
  sur\-faces in  superspace.
  Theor. Mat. Fiz. {\bf 52}, (1982), 375--383.
\m
[10] Baranov,M.A., Schwarz,A.S.:  Characteristic Clases of
Supergauge Fields {\it Funkts. Analiz i ego pril.},
 {\bf 18} No.2, (1984), 53--54.
 \smallskip
 Cohomologies of supermanifolds.
 {\it Funkts. analiz i ego pril.}.
 {\bf 18} No.3, (1984), 69--70.
\m
 [11] Voronov, T.:  Geometric Integration Theory on supermanifolds.
   {\it Sov.Sci.Rev.C Math.} {\bf 9}, (1992), 1--138.
      \m
 [12] A.S.Schwarz-  Supergravity, complex geometry and
 G-structures
   Commun. Math. Phys., {\bf 87}, (1982) 37--63 .
 \m
[13] Shander, V.N.:
  Analogues of the Frobenius and Darboux Theorems for Supermanifolds.
  {\it Comptes rendus de l'
 Academie bulgare des Sciences}, {\bf 36}, n.3,
 (1983), 309--311.
  \m
 [14] Dubrovin,B.A., Novikov,S.P., Fomenko,A.T.:
 {\it  Modern Geometry. Methods and Applications}, Moscow, Nauka, 1986.
\bye